\documentclass [11pt]{amsart}
\usepackage{epic,eepic,latexsym, amssymb, amscd, amsfonts}

\input xy
\xyoption {all}

\def \nquot {{{\text{Quot}}_{X} ({\mathcal O}^{k+r}, r, d)}}
\def \jac {{\text{Jac}}}
\def \m {{\mathfrak M}}

\newtheorem{theorem}{Theorem}

\newtheorem{conjecture}{Conjecture}
\newtheorem {corollary}{Corollary}
\newtheorem {proposition}{Proposition}

\theoremstyle{definition}

\newtheorem {example}{Example}

\theoremstyle {definition}

\begin{document}

\title[A tour of theta dualities on moduli spaces of sheaves]{A tour of theta dualities on moduli spaces of sheaves}
\author {Alina Marian}
\address {School of Mathematics}
\address {Institute for Advanced Study}
\email {marian@math.ias.edu}
\author {Dragos Oprea}
\address {Department of Mathematics}
\address {Stanford University}
\email {oprea@math.stanford.edu}
\date{}
\maketitle

\begin {abstract} The purpose of this paper is twofold. First, we survey known results about theta dualities on moduli spaces of sheaves on
curves and surfaces. Secondly, we establish new such dualities in the surface case. Among others, the case of elliptic K3
surfaces is studied in detail; we propose further conjectures which are shown to imply strange duality. \end {abstract}

\section{Introduction}

The idea that sections of the determinant line bundles on moduli spaces of sheaves
are subject to natural dualities was first formulated, and almost exclusively pursued, in the case 
of vector bundles on curves. The problem can however be generally stated as follows.

Let $X$ be a smooth complex projective curve or surface with polarization $H$, let $v$ be a class in the Grothendieck group $K(X)$ of coherent 
sheaves on $X$. Somewhat imprecisely, we denote by $\m_v$ the moduli space of Gieseker $H$-semistable sheaves on $X$ of class $v$. Consider the 
bilinear form in $K$-theory \begin{equation} \label{euler} (v, w) \mapsto \chi (v \otimes w), \, \, \text{for}\, \, v, w \in K(X), \end{equation} and let 
$v^{\perp} \subset K(X)$ consist of the $K$-classes orthogonal to $v$ relative to this form.  There is a group homomorphism $$\Theta: v^{\perp} 
\longrightarrow \text{Pic} \, \m_v, \, \, \, w \mapsto \Theta_w,$$ which was extensively considered in \cite {dn} when $X$ is a curve, and also in 
\cite {lp}, \cite{jun} in the case when $X$ is a surface, as part of the authors' study of the Picard group of the moduli spaces of sheaves. When $\m_v$ admits a universal sheaf ${\mathcal E} \rightarrow \m_v \times X$,  $\Theta_w$ is given by $$\Theta_w = \det {\bf R} p_{!}(\mathcal E \otimes q^{\star} F)^{-1}.$$ Here $F$ is any sheaf with $K$-type $w$, and $p$ and $q$ are the two projections from $\m_v \times X$. The line bundle $\Theta_w$ is well defined also when the moduli space $\m_v$ is not fine, by descent from the Quot scheme.

Consider now two classes $v$ and $w$ in $K(X)$, orthogonal with respect to the bilinear form \eqref{euler} {\it i.e.,} satisfying $\chi (v \otimes w) 
= 0.$ Assume that for any points $[E] \in {\mathfrak M}_v$ and $[F] \in {\mathfrak M}_w$, the vanishing, vacuous when $X$ is a curve, \begin 
{equation}\label{vanish}H^2 (E \otimes F) = 0\end {equation} occurs. Suppose further that the locus \begin{equation} \label{theta} \Theta = \{ (E, F) 
\in {\mathfrak M}_v \times {\mathfrak M}_w \, \, \text{with} \, \, h^0 (E \otimes F) \neq 0 \} \end{equation} gives rise to a divisor of the line bundle 
$\mathcal L$ which splits as \begin{equation}\label{split}{\mathcal L} = \Theta_w \boxtimes \Theta_v \, \, \, \text{on} \,\, \m_v \times 
\m_w.\end {equation} We then obtain a
morphism, well-defined up to scalars, \begin{equation} \label{duality} {\mathsf D}: {H^0 (\m_v, \Theta_w)}^{\vee} \longrightarrow H^0 (\m_w,
\Theta_v). \end{equation} \vskip.1in

The main questions of geometric duality in this context are simple to state and fundamentally naive.\vskip.1in

{\it {\bf Question 1.} What are the constraints on $X$, $v,$ and $w$ subject to which one has, possibly with suitable 
variations in the meaning of $\m_v$, 
\begin{equation}
\label{equality}
h^0 (\m_v, \Theta_w) = h^0 (\m_w, \Theta_v)?
\end{equation}

{\bf Question 2.} In the cases when the above equality holds, is the map ${\mathsf D}$ of equation \eqref{duality} an 
isomorphism? } \vskip.1in

This paper has two goals. One is to survey succinctly the existent results addressing Questions $1$ and $2.$ The other one is to study the two
questions in new geometric contexts, providing positive answers in some cases and support for further conjectures in other cases. Throughout we
refer to the isomorphisms induced by jumping divisors of type \eqref{theta} as {\it theta dualities}, or {\it strange dualities}. The latter term is in keeping
with the terminology customary in the context of moduli spaces of bundles on curves.

We begin in fact by reviewing briefly the arguments that establish the duality in the case of vector bundles on curves. We point out that in this
case the isomorphism can be regarded as a generalization of the classical Wirtinger duality on spaces of level $2$ theta functions on abelian
varieties. We moreover give a few low-rank/low-level examples. We end the section devoted to curves by touching on Beauville's proposal \cite{Bea2} concerning a strange duality on moduli spaces of symplectic bundles; we illustrate the symplectic duality by an example.

The rest of the paper deals with strange dualities for moduli spaces of sheaves on surfaces. In this context, Questions $1$ and $2$ were first
posed by Le Potier. Note that in the curve case there are strong representation-theoretic reasons to expect affirmative 
answers to these questions. By contrast, no analogous reasons are known to us in the case of surfaces.

The first examples of theta dualities on surfaces are given in Section 3. There, we explain the theta isomorphism for pairs of rank 1 moduli
spaces {\it i.e.,} for Hilbert schemes of points. We also give an example of strange duality for certain pairs of rank 0 moduli spaces.

Section 4 takes up the case of moduli spaces of sheaves on surfaces with trivial canonical bundle. The equality \eqref{equality} of dimensions for
dual spaces of sections has been noted \cite{ogr2}\cite {gny} in the case of sheaves on $K3$ surfaces. Moreover, it was recently established in a few different flavors for sheaves on abelian surfaces \cite{mo3}. We review the numerical statements, which lead one to speculate that the duality map is an
isomorphism in this context. We give a few known examples on $K3$ surfaces, involving cases when one of the two moduli spaces is itself a $K3$
surface or a Hilbert scheme of points. It is likely that more general instances of strange duality can be obtained, starting with the
isomorphism on Hilbert schemes presented in Section 3 and applying Fourier-Mukai transformations. We hope to investigate this elsewhere. Finally,
the last part of this section leaves the context of trivial canonical bundle, and surveys the known cases of strange duality on the projective plane,
due to D\v{a}nil\v{a} \cite {D1} \cite{D}.

Section 5 is devoted to moduli spaces of sheaves on elliptically fibered K3 surfaces. We look at the case when the first Chern class of the sheaves
in the moduli space has intersection number 1 with the class of the elliptic fiber. These moduli spaces have been explicitly shown birational to
the Hilbert schemes of points on the same surface \cite{ogr}. We conjecture that strange duality holds for many pairs of such moduli spaces,
consisting of sheaves of ranks at least two. We support the conjecture by proving that if it is true for one pair of ranks, then it holds for any
other, provided the sum of the ranks stays constant. Since any moduli space of sheaves on a $K3$ surface is deformation equivalent to a moduli space on an elliptic $K3$, the conjecture has
implications for generic strange duality statements.

\subsection {Acknowledgements} We would like to thank Jun Li for numerous conversations related to moduli spaces of sheaves and Mihnea Popa for bringing to our attention the question of strange dualities on surfaces. We are grateful to Kieran O'Grady for clarifying a technical point. Additionally, we acknowledge the financial support of the NSF. A.M. is very grateful to Jun Li and the Stanford Mathematics Department for making possible a great stay at Stanford in the spring of 2007, when this article was started. 

\section{Strange duality on curves}

\subsection {General setup} When $X$ is a curve, the topological type of a vector bundle is given by its rank and degree, and its class in the Grothendieck group by its rank
and determinant. We let $\m (r,d)$ be the moduli space of semistable vector bundles with fixed numerical data given by the rank $r$ and degree $d$.
Similarly, we denote by $\m (r,\Lambda)$ the moduli space of semistable bundles with rank $r$ and fixed
determinant $\Lambda$ of degree $d$. For a vector bundle $F$ which is orthogonal to the bundles in $\m(r,d)$ {\it i.e.,} $$\chi(E\otimes F)=0,
\text { for all } E\in \m(r,d),$$ we consider the jumping locus $$\Theta_{r,F}=\{E\in \m(r,d) \text { such that } h^{0}(E\otimes F)\neq 0\}.$$ The
additional numerical subscripts of thetas indicate the ranks of the bundles that make up the corresponding moduli space.

The well-understood structure of the Picard group of $\m (r,d)$, given in \cite {dn}, reveals immediately that the line bundle associated to
${\Theta}_{r, F}$ depends only on the rank and determinant of $F$. Moreover, on the moduli space $\m(r, \Lambda)$, $\Theta_{r,F}$ depends only on the
rank and degree of $F$. In fact, the Picard group of $\m(r, \Lambda)$ has a unique ample generator $\theta_{r}$, and therefore $${\Theta}_{r,F}\cong \theta_r^l,$$ for some integer $l$.

For numerical choices $(r,d)$ and$(k,e)$ orthogonal to each other, the construction outlined in the Introduction gives a
duality map \begin {equation}\label{sdm1}{\mathsf {D}}: {H^0 (\m (r,\Lambda), {\Theta}_{r,F} )}^{\vee} \longrightarrow H^0 (\m (k,e), {\Theta}_{k, E}),
\, \, \text{for} \, \, E \in \m (r,\Lambda), F \in \m (k,e).\end {equation} The {\it strange duality} conjecture of Beauville \cite{bsurvey}
and Donagi-Tu \cite{dt} predicted that the morphism $\mathsf {D}$ is an isomorphism. This was recently proved for a generic smooth curve in \cite{bel1}, which inspired a subsequent argument  of \cite{mo1} for all smooth curves. The statement for all curves also follows from the generic-curve case in conjunction with the recent results of \cite{bel2}. We briefly review the arguments in the subsections below.

\subsection {Degree zero} The duality is most simply formulated when $d = 0$, $e = k (g-1)$ on the moduli space $\m (r,\mathcal O)$ of rank $r$ bundles 
with trivial determinant. There is a canonical line bundle ${\Theta_k}$ on the moduli space $\m (k, k(g-1)),$ associated 
with the divisor $$\Theta_{k}=\{F \in \m (k, k(g-1)) \, \, \text{such that} \, \, h^0 (F) = h^1 (F) \neq 0\}.$$ 
The map \eqref{sdm1} becomes 
\begin{equation}
\label{sd}
{\mathsf {D}}: H^0 (\m (r,\mathcal O), \theta_r^k )^{\vee} \longrightarrow H^0 (\m (k, k(g-1)), \Theta_k^r),
\end{equation} 
Note that this interchanges the rank of the bundles that make up the moduli space and the level (tensor power) of the 
determinant line bundle on the moduli space.\vskip.1in

To begin our outline of the arguments, let 
\begin{equation}
v_{r,k} = \chi (\m (r,\mathcal O), \theta_r^k )
\end{equation}
be the {\it {Verlinde number}} of rank $r$ and level $k$. As the theta 
bundles have no higher cohomology, the Verlinde number computes in
fact the dimension $h^{0}(\theta_{r}^{k})$ of the space of sections. The most elementary formula for $v_{r,k}$ reads

\begin{equation} \label{elem} v_{r,k} = \frac{r^{g}}{(r+k)^{g}}\sum_{\stackrel{S\sqcup T=\{1, \ldots, r+k\}}{|S|=k}} \prod_{\stackrel{s\in S}{t\in T}}
\left|2\sin\pi\frac{s-t}{r+k}\right|^{g-1}. \end{equation} This expression for $v_{r,k}$ was established through a concerted effort and variety
of approaches that spanned almost a decade of work in the moduli theory of bundles on a curve. It is beyond the purpose of this note to give an overview of
the Verlinde formula. Nonetheless, let us mention here that \eqref{elem} implies the symmetry $$h^{0}(\m(r, \mathcal O), \Theta_{r}^{k})= h^0 (\m
(k, k(g-1)), \Theta_k^r),$$ required by Question $1$, setting the stage for the strange duality conjecture. 
\vskip.1in \begin {example}
{\it Level 1 duality.} Equation \eqref{elem} simplifies dramatically in level $k=1,$ yielding $$v_{r,1}=\frac{r^{g}}{(r+1)^{g-1}} \prod_{p=1}^{r}\left|2 \sin
\frac{p\pi}{r+1}\right|^{g-1}=r^{g}.$$ This coincides with the dimension of the space of level $r$ classical theta functions on the Jacobian. The
corresponding isomorphism \begin{equation}\label{thh}\mathsf D:H^{0}(\m(r, \mathcal O), \theta_{r})^{\vee}\to H^{0}(\jac^{g-1}(X),
\Theta^{r}_{1})\end{equation} was originally established in \cite {bnr}, at that time in the absence of the Verlinde formula. Nonetheless, once the Verlinde formula is known, the isomorphism \eqref{thh} may be proved by an easy argument which we learned from Mihnea Popa; see also \cite {bel2}. Let $\jac[r]$ denote the group of $r$-torsion points on the Jacobian, and let $\mathsf H[r]$ be the Heisenberg group of
the line bundle $\Theta_{1}^{r}$, exhibited as a central extension $$0\to \mathbb C^{\star}\to \mathsf H[r]\to \jac[r]\to 0.$$ Now $H^{0}(\jac^{g-1}(X),
\Theta_{1}^{r})$ is the Schr\"odinger representation of the Heisenberg group $\mathsf H[r]$ {\it i.e.,} the unique irreducible representation on which the 
center
$\mathbb C^{\star}$ acts by homotheties. We argue that the left hand side of \eqref{thh} also carries such a representation of $\mathsf H[r]$. Indeed, the
tensor product morphism $$t:\m(r, \mathcal O)\times \jac^{g-1}(X)\to \m(r, r(g-1))$$ is invariant under the action of $\jac[r]$ on the source, given by
$$(\zeta, (E, L))\mapsto (\zeta^{-1}\otimes E, \zeta\otimes L), \, \, \text{for} \, \zeta \in \jac[r], \, (E, L) \in \m (r, \mathcal O) \times \jac^{g-1} 
(X).$$
The pullback  $$t^{\star}\Theta_{r} =\theta_{r}\boxtimes \Theta^{r}_{1}$$ carries an action of $\jac[r]$, while the theta bundle $\Theta_{1}^{r}$ on
the right carries an action of $\mathsf H[r]$. This gives an induced action of the Heisenberg group on $\theta_{r},$ covering the tensoring action of
$\jac[r]$ on $\m(r, \mathcal O)$. The non-zero morphism $\mathsf D$ is clearly $\mathsf H[r]$-equivariant, hence it is an isomorphism.\end{example}

\vskip.1in

It is more difficult to establish the isomorphism \eqref {sd} for arbitrary levels. The argument in \cite {bel1} starts by noticing
that since the dimensions of the two spaces of sections involved are the same, and since the map is induced by the divisor \eqref{theta}, ${\mathsf
{D}}$ is an isomorphism if one can generate pairs 
\begin{equation}
\label{pairs0}
(E_i, F_i) \in \m (r,\mathcal O) \times \m (k, k(g-1)), \, \, \, 1\leq i \leq v_{r,k} 
\end{equation}
such that $$E_{i} \in \Theta_{F_{j}} \text { if and only if } i\neq j.$$ Equivalently, one requires
\begin{equation}
\label{req}
h^0 (E_i \otimes F_j) = 0 \,\, \text{if and only if} \, \, i = j.
\end{equation} Finding the right number of such pairs relies on giving a suitable enumerative interpretation to the Verlinde 
formula \eqref{elem}, and this is the next step in \cite{bel1}. 
The requisite vector bundles \eqref{pairs0} satisfying \eqref{req} are assembled first on a rational nodal curve of genus $g$. They are then
deformed along with the curve, maintaining the same feature on neighboring smooth curves, and therefore establishing the duality for generic
curves. To obtain the right count of bundles on the nodal curve, the author draws benefit from the {\it fusion rules} \cite{ueno} of the Wess-Zumino-Witten theory, which express the Verlinde numbers in genus $g$ in terms of Verlinde numbers in lower genera, and with point insertions. The latter are Euler characteristics of line bundles on moduli spaces of bundles with parabolic structures at the given points.  Belkale realizes the
count of vector bundles with property \eqref{req} on the nodal curve as an enumerative intersection in Grassmannians with recursive traits
precisely matching those of the fusion rules.

Recently, the author extended the result from the generic curve case to that of arbitrary curves \cite {bel2}. He considers the spaces of
sections $H^0 (\m (r, \mathcal O), {\mathcal L}^k )$ relatively over families of smooth curves. These spaces of sections give rise to vector bundles which come equipped with
Hitchin's projectively flat connection \cite {hitchin}. Belkale proves, importantly, that \begin {theorem}\cite {bel2} The relative strange duality
map is projectively flat with respect to the Hitchin connection. In particular, its rank is locally constant. \end {theorem} This result moreover
raises the question of extending the strange duality isomorphism to the boundary of the moduli space of curves. \vskip.1in

The alternate proof of the duality in \cite{mo1} is inspired by \cite{bel1}, in particular by the interpretation of the 
isomorphism as solving a counting problem for bundles satisfying  \eqref{pairs0},  \eqref{req}.  This count is carried out in \cite{mo1} using the close intersection-theoretic rapport between the
moduli space of bundles $\m (r,d)$ and the Grothendieck Quot scheme $\text{Quot}_{X} ({\mathcal O}^{N}, r, d)$ parametrizing rank $r$, degree $d$
coherent sheaf quotients of the trivial sheaf ${\mathcal O}^N$ on $X$. The latter is irreducible for large degree $d$, and compactifies the space
of maps $\text{Mor}_{d}(X, G(r, N))$ from $X$ to the Grassmannian of rank $N-r$ planes in ${\mathbb C}^N$. The $N$-asymptotics of tautological
intersections on the Quot scheme encode the tautological intersection numbers on $\m (r,d)$ \cite{mo2}. This fact prompts one to take the viewpoint
that the space of maps from the curve to the classifying space of $GL_{r}$ is the right place to carry out the intersection theory of the moduli
space of $GL_{r}$-bundles on the curve.

As a Riemann-Roch intersection on $\m (r, \mathcal O)$, the Verlinde number $v_{r,k}$ 
is also expressible on the Quot scheme, most 
simply as a tautological intersection on $\nquot$ which parametrizes rank $r$ quotients of ${\mathcal O}^{k+r}$ of suitably large degree $d$. Precisely, let 
${\mathcal E}$ denote the rank $k$ universal subsheaf of ${\mathcal O}^{k+r}$ on $\nquot \times X,$
and set $$a_k = c_k ({{\mathcal E}^{\vee} |}_{\nquot \times \{\text{point} \}}).$$ Then

\begin{equation}
\label{quot}
v_{r,k} = h^0 (\m (r, \mathcal O), \theta_r^k) = \frac{r^g}{(k+r)^g} \int_{\nquot} a_k^{\text{top}}.
\end{equation} 

This formula was essentially written down by Witten \cite{witten} using physical arguments; it reflects the relationship between the $GL_{r}$ WZW
model and the sigma model of the Grassmannian, which he explores in \cite{witten}. Therefore, the modified Verlinde number
$${\tilde{v}}_{r, k} = \frac{(k+r)^{g}}{r^{g}}v_{r,k}=\int_{\nquot} a_k^{\text{top}}$$ becomes a count of Quot scheme points which obey incidence
constraints imposed by the self intersections of the tautological class $a_k$. Geometrically interpreted, these constraints single out the finitely
many subsheaves $$E_i \hookrightarrow {\mathcal O}^{k+r}$$ which factor through a subsheaf $S$ of ${\mathcal O}^{k+r}$ of the same rank $k+r$ but
of lower degree. Therefore, one obtains diagrams \begin{equation} \label{zeroquot} {\xymatrix{ 0 \ar[r] & E_i \ar[r] \ar[dr] & S \ar [r] \ar [d] &
F_i \ar[r] & 0 \\ & & {\mathcal O}^{k+r} & & }}, \, \, \, \text{for} \, \, 1\leq i \leq {\tilde{v}}_{r,k}, \end{equation} where the top sequence is
exact and the triangle commutes. The sheaf $S$ is obtained by a succession of elementary modifications of ${\mathcal O}^{k+r}$ dictated by choices
of representatives for the class $a_k$. For generic choices, the $\tilde v_{r,k}$ exact sequences \eqref{zeroquot} are the points of a {\it{smooth}}
zero-dimensional Quot scheme on $X$ associated with $S$, and therefore automatically satisfy $$h^0 (E_i^{\vee} \otimes F_j) = 0 \, \, \, \text{iff}
\, \, \, i = j.$$ 

Now the modified Verlinde number $\tilde{v}_{r,k}$ which gives the number of such pairs is the Euler characteristic of a
{\it{twisted}} theta bundle on $\m (r,0).$ Precisely, for any reference line bundles $L$ and $M$ on $X$ of degree $g-1$, $$\tilde{v}_{r,k} = h^0
(\m (r, 0), \Theta_{r, M}^k \otimes {\det}^{\star} \Theta_{1,L} ), $$ where $$\det: \m (r,0)\to \jac^0(X)$$ is the morphism which takes bundles to their determinants.

Along the lines of the remarks preceding equation \eqref{req}, the pairs of bundles $(E_i, F_i)$ coming from 
\eqref{zeroquot} easily give the following result:

\begin {theorem} \label{mainsd} \cite {mo1} (Generalized Wirtinger duality.) For any line bundles $L$ and $M$ of degree $g-1$, there is an isomorphism 
\begin{equation}
\label{sd1alt}
\widetilde {\mathsf {D}}:H^{0}\left (\m(r,0),   
\Theta_{r,M}^{k}\otimes {\det}^{\star} \Theta_{1,L}\right)^{\vee} \to H^{0} \left (\m (k,0), \Theta_{k,M}^{r}\otimes
\left ((-1)\circ {\det}\right)^{\star}
\Theta_{1,L} \right ).
\end{equation}
Here, $(-1)$ denotes the multiplication by $-1$ on the Jacobian $\jac^{0}(X)$.
\end {theorem}

\begin {example} {\it Wirtinger duality.} When $k = r= 1$, the twisted duality map $\widetilde {\mathsf D}$ is the
classical Wirtinger duality on the abelian variety $\jac(X)$ \cite {mumford}. 
For simplicity let us take $L = M$, and let us assume that $L$ is a theta characteristic. Then the line bundle
$$\Theta = \Theta_{1, L} = \Theta_{1, M}$$ is symmetric {\it i.e.,} $(-1)^{\star} \Theta \cong \Theta.$ Let $\delta$ be the morphism $$\jac(X)\times 
\jac(X)\ni(A,
B)\to (A^{-1}\otimes B, A\otimes B)\in \jac(X)\times \jac(X).$$ 
The see-saw theorem shows that
$$\delta^{\star}\left(\Theta\boxtimes \Theta\right)=2\Theta \boxtimes 2\Theta.$$ In turn, this induces the Wirtinger self-duality $$\mathsf
W=\widetilde {\mathsf D}:\left |2\Theta\right|^{\vee}\to \left|2\Theta \right|.$$ One may establish the fact that $\mathsf W$ is an isomorphism by
considering as before the action of the Heisenberg group $\mathsf H[2]$. We refer the reader to Mumford's article \cite {mumford} for such an argument. Note also that the study of the duality morphism $\mathsf W$ leads to various addition formulas for theta functions;
these have consequences for the geometric understanding of the $2\Theta$-morphism mapping the abelian variety into projective space. \end
{example}

\begin {example} {\it Rank 2 self-duality.} The next case of self-duality occurs in rank $2$ and level $2$. The morphism $\widetilde {\mathsf {D}}$ leads to
the study of the line bundles $$\Theta_{L, M}=\Theta_{2,M}^{2}\otimes {\det}^{\star}\Theta_{1,L},$$ where $L$ and $M$ are line bundles of degree $g-1$ on
$X$. We assume as before that $\Theta_{1,L}$ is symmetric. We claim that the line bundles $\Theta_{L,M}$ are globally generated, giving 
rise to morphisms $$f_{L, M}:\m(2, 0)\to |\Theta_{L, M}|^{\vee}\cong |\Theta_{L, M}|.$$ It would be
interesting to understand the geometry of this self-duality in more detail.

To see that $\Theta_{L,M}$ is globally generated, observe that if $A_{1}, A_{2}$ are degree $0$ line bundles on $X$, then, using the formulas in
\cite {dn}, we have $$\Theta_{L,M}=\Theta_{2,M\otimes A_{1}}\otimes \Theta_{2, M\otimes A_{2}}\otimes {\det}^{\star}\Theta_{1, L\otimes
A_{1}^{\vee}\otimes A_{2}^{\vee}}.$$ This rewriting of the line bundle $\Theta_{L,M}$ gives a section vanishing on the locus $$\Delta_{A_{1},
A_{2}}=\left\{E\in \m (2,0) \text { such that } h^0(E\otimes M\otimes A_{1})\neq 0 \text { or } h^{0}(E\otimes M\otimes A_{2})\neq 0 \text { or
}\right.$$ $$\left. h^{0}(\det E\otimes L\otimes A_{1}^{\vee}\otimes A_{2}^{\vee})\neq 0\right\}.$$ Fixing $E\in \m (2,0)$, it suffices to show that one
of the sections $\Delta_{A_{1}, A_{2}}$ does not vanish at $E$. This amounts to finding suitable line bundles $A_{1}, A_{2}$ such that $E\not
\in\Delta_{A_{1}, A_{2}},$ or equivalently $$A_{1}, A_{2}\not \in \Theta_{1, E\otimes M} \text { and } A_{1}^{\vee}\otimes A_{2}^{\vee}\not \in
\Theta_{1, \det E\otimes L}.$$ The existence of $A_{1}, A_{2}$ is guaranteed by Raynaud's observation \cite{R} that for any vector bundle $V$ of
rank $2$ and degree $2g-2$, the locus $$\Theta_{1, V}=\left\{A \in \text {Jac} (X), h^{0}(V\otimes A)=h^1(V\otimes A)\neq 0\right\}$$ has
codimension $1$ in $\text {Jac}(X)$.  \end {example}

Turning to the original formulation \eqref{sd} of strange duality on curves, Theorem \ref{mainsd} and a restriction argument on the Jacobian imply that  
\begin
{theorem} \cite {mo1} The strange duality map $\mathsf {D}$ of \eqref{sd} is an isomorphism.  \end {theorem} As a consequence, we obtain the following 

\begin {corollary} \label{cor1} The theta divisors $\Theta_{E}$ generate the linear system $\left|\theta_{r}^{k}\right|$ on the moduli space
$\mathfrak M(r, \mathcal O)$, as $E$ varies in the moduli space $\mathfrak M(k, k(g-1))$.  \end {corollary}

\begin {example} {\it Verlinde bundles and Fourier-Mukai.} As observed by Popa \cite {po}, the strange duality isomorphism \eqref{sd} can be elegantly
interpreted by means of the Fourier-Mukai transform. Fixing a reference theta characteristic $L$, Popa defined the Verlinde bundles $$\mathsf
E_{r, k}={\det}_{\star}\left(\Theta_{r, L}^{k}\right)$$ obtained as the pushforwards of the $k$-pluritheta bundles via the morphism $$\det:\m(r, 0)\to
\jac(X).$$ He further showed that the Verlinde bundles satisfy $IT_{0}$ with respect to the normalized Poincar\'{e} bundle on the Jacobian, and studied
the Fourier-Mukai transform $\widehat {\mathsf E}_{r,k}$. The construction of the duality morphism $\mathsf D$ globalizes as we let the determinant
of the bundles in the moduli space vary in the Jacobian. The ensuing morphism $$\mathsf {D}: \mathsf E_{r,k}^{\vee}\to \widehat {\mathsf E}_{k,r}$$
collects all strange duality morphisms for various determinants, and it is therefore an isomorphism. We moreover note here that in level $k=1$,
Popa showed that both bundles are stable, therefore giving another proof of the fact that $\mathsf {D}$ is an isomorphism.

\end {example} \subsection{Arbitrary degrees} The proof outlined above effortlessly generalizes to arbitrary degrees. Specifically, let $r$ and $d$ be coprime
integers, and $h$, $k$ be any two non-negative integers, and fix $S$ an auxiliary line bundle of degree $d$ and rank $r$.

\begin{theorem}\cite{mo1} \label{arbdegree} There is a level-rank duality isomorphism between $$H^{0}\left(\m(hr, hd), \Theta_{hr, M\otimes
S^{\vee}}^{k}\otimes {\det}^{\star} \Theta_{1, L\otimes (\det S)^{-h}}\right) \text { and }$$ $$H^{0}\left(\m(kr, -kd), \Theta_{kr, M\otimes S}^{h}\otimes {((-1)
\circ \det)}^{\star}\Theta_{1, L\otimes \left(\det S\right)^{-k}}\right).$$ \end{theorem}

Finally, just as in \eqref{sd}, the tensor product map induces the strange duality morphism for arbitrary degree:
\begin{equation}\label{sdd}\mathsf{D}:H^{0}\left(\m\left(hr, (\det S)^{h}\right), \theta^k_{hr}\right)^{\vee}\to
H^{0}\left(\m\left(kr, k(r (g-1) - d)\right), \Theta_{kr,S}^{h}\right). \nonumber\end{equation} 

\begin{theorem} \cite {mo1} The strange duality morphism $\mathsf {D}$ is an isomorphism. 
\end {theorem} 

\subsection {Symplectic strange duality} It is natural to inquire if the same duality occurs for moduli spaces of principal bundles with arbitrary structure groups. The symplectic group should be considered next, due to the fact that the moduli space of symplectic bundles is locally factorial.
A recent conjecture of Beauville \cite {Bea2}, which we now explain, focuses on this case. 

Consider the moduli space $\mathcal M_{Sp_r}$ of
pairs $(E, \phi)$, where $E$ is a semistable bundle of rank $2r$, such that $$\det E=\mathcal O_{X} \text { and } \phi:\Lambda^{2} E\to \mathcal O_{X}
\text{ is a non-degenerate alternate form}.$$ Similarly, let $\widehat{\mathcal M}_{Sp_k}$ denote the cousin moduli space of pairs $(F, \psi)$ as
above, with the modified requirements $$\det F=K_{X}^k, \text { and } \psi: \Lambda^{2} F\to K_X \text { is a nondegenerate alternate form}.$$ We
let $\mathcal L_r$ and $\widehat {\mathcal L}_{k}$ be the determinant bundles on ${\mathcal {M}}_{Sp_{r}}$ and $\widehat{\mathcal M}_{Sp_k}$ respectively. Now observe the tensor product map $$t:\mathcal
M_{Sp_r}\times \widehat {\mathcal M}_{Sp_k}\to \widehat{\mathcal M}^{+}$$ given by $$(E, \phi) \times (F, \psi)\to (E\otimes F, \phi\otimes
\psi).$$ Its image is contained in the moduli space $\widehat{\mathcal M}^{+}$ of even orthogonal pairs $(G,q),$ consisting of semistable bundles
$G$ of rank $4rk$, determinant $K_{X}^{2rk}$, with $h^{0}(G)$ even, and endowed with a quadratic form $$q:\text {Sym}^{2} G\to K_{X}.$$

Beauville showed that the pullback divisor $t^{\star}\Delta$, with $$\Delta= \left\{(G,q)\in \widehat{\mathcal M}^{+} \text { such that }
h^{0}(G)\neq 0\right\},$$ determines a canonical section of $\mathcal L_r^{k}\boxtimes \widehat {\mathcal L}_{k}^{r}$. Hence, it gives rise to a
morphism $$\mathsf {D}: H^{0}(\mathcal M_{Sp_r}, \mathcal L_r^{k})^{\vee}\to H^{0}(\widehat{\mathcal M}_{Sp_{k}},\widehat{\mathcal L}_k^{r}).$$
Beauville checked that the dimensions of the spaces of sections match, $$h^{0}(\mathcal M_{Sp_r}, \mathcal L_r^{k})=h^{0}(\widehat{\mathcal
M}_{Sp_{k}},\widehat{\mathcal L}_k^{r}),$$ and further conjectured that

\begin {conjecture} \cite {Bea2} The morphism $\mathsf {D}$ is an isomorphism.
\end {conjecture}

\begin {example} {\it Level $1$ symplectic duality}. The particular case $k=1$ of this conjecture is a consequence of the usual strange
duality theorem. In this situation, there is an isomorphism $$(\widehat {\mathcal M}_{Sp_{1}}, \widehat {\mathcal L}_{1})\cong (\mathfrak M(2, K_{X}),
\theta_{2}).$$ It suffices to explain that the linear system $|\theta^{r}_{2}|$ is spanned by the sections $$\Delta_{E}=\{(F, \psi) \text { such
that } h^{0}(E\otimes F)\neq 0\},$$ as $E$ varies in $\mathcal M_{Sp_{r}}.$ In fact, we only need those $E$'s of the form $$E=E'\oplus E'^{\vee},$$
for $E'\in \mathfrak M(r, 0).$ In this case, Beauville observed that $$\Delta_{E}=\Theta_{E'}=\{F\in \m (2, K_X), \text { such that } h^{0}(F\otimes
E')\neq 0 \}.$$ These generate the linear system $|\theta_{2}^{r}|$, by Corollary \ref{cor1}.\end {example}

\section{Duality on Hilbert schemes of points on a surface} \label{1-1} In this section, we start investigating strange duality phenomena for
surfaces. We begin with two examples involving the Hilbert scheme of points. In the next sections we will use these basic cases to obtain new examples of theta dualities for surfaces with trivial canonical bundles.

\subsection{Notation.}
\label{hilbert1}

Denote by $X^{[k]}$ the Hilbert scheme of $k$ points on the projective surface $X$, and let
$$ 0 \rightarrow {\mathcal I}_{\mathcal Z} \rightarrow {\mathcal O} \rightarrow {\mathcal O}_{\mathcal Z} 
\rightarrow 0 $$ be the universal family on $X^{[k]} \times X.$ Let $p$ be the projection from the product 
$X^{[k]} \times X$ to the Hilbert scheme, and $q$ be the projection to the surface $X$.

For any sheaf $F$ on $X$, let $F^{[k]}$ be the line bundle $$F^{[k]} = \det p_{!} ({\mathcal O}_{\mathcal Z} \otimes q^{\star} F)  =
\left(\det (p_{!} ({\mathcal I}_{\mathcal Z} \otimes q^{\star} F) ) \right )^{-1}.$$ (Note that in the literature on Hilbert schemes of points,
$F^{[k]}$ often refers to the pushforward itself, not its determinant. It will be convenient for us not to comply with this practice, and single
out, as above, the {\it{determinant}} line bundle by this notation.) Further, for any line bundle $N$ on $X$, the $S_k$-equivariant line bundle
$N^{\boxtimes k}$ on $X^k$ descends to a line bundle on the symmetric product $X^{(k)}$. Let $N_{(k)}$ be the pullback of this descent bundle on
the Hilbert scheme $X^{[k]}$, under the Hilbert-Chow morphism $$f: X^{[k]} \rightarrow X^{(k)}.$$
We recall from \cite {EGL} that
\begin{equation}
\label{picard}
F^{[k]} = (\det F)_{(k)} \otimes M^{\, \text{rk}\, F}, 
\end{equation}
where we set
$$M = {\mathcal O}^{[k]}.$$
Equation \eqref{picard} implies in particular that
\begin{equation}
\label{seesaw}
F^{[k]} \cong (F \otimes I_Z)^{[k]}\end {equation} for any zero-cycle $Z$, since $F^{[k]}$ only depends on the determinant and rank of $F$. 

\subsection{The strange duality isomorphism.}
\label{hilbert2}

To set up the strange duality map, let $L$ be a line bundle on $X$ without higher cohomology. Let $$n = \chi (L) = h^0 (L),$$ and pick $1\leq k\leq n.$ Denote by $\theta_{L,k}$ the
divisor on $X^{[k]} \times X^{[n-k]}$ which, away from the codimension $2$ locus of pairs $(Z, W)$ with overlapping support,  is given by $$\theta_{L,k} = \left \{ (I_Z, I_W) \in X^{[k]} \times X^{[n-k]} \, \text{such that} \, \,h^0 (I_Z
\otimes I_W \otimes L) \neq 0 \right \}.$$
Using the seesaw theorem and equation \eqref{seesaw}, we find that 
$${\mathcal O} (\theta_{L,k}) \cong L^{[k]} \otimes L^{[n-k]} \, \, \text{on} \, \, X^{[k]} \times
X^{[n-k]}.$$
In particular,
$$\theta_{L,k} \in H^0 (X^{[k]}, L^{[k]} ) \otimes H^0 (X^{[n-k]}, L^{[n-k]}).$$
In this subsection, we record the following

\begin{proposition}\label{hilbert} Assume that $L$ is a line bundle with $\chi(L)=n \geq k$, and no higher cohomology. The map $${\mathsf D}_{L}: 
H^0
(X^{[k]}, L^{[k]})^{\vee} \longrightarrow H^0 (X^{[n-k]}, L^{[n-k]})$$ induced by the divisor $\theta_{L,k}$ is an isomorphism. \end{proposition}

{\it{Proof.}} The space of sections $H^0 (X^{[k]}, L^{[k]})$ can be realized explicitly in terms of sections of $L$ on $X$, cf. \cite {EGL}. The Proposition follows from this identification. Specifically $H^0 (X^{[k]}, L^{[k]})$ can be viewed as the invariant part of $H^0 (X,
L)^{\otimes k}$ under the antisymmetric action $\epsilon$ of the permutation group $S_k$, \begin{equation} \label{identify} H^0 (X^{[k]}, L^{[k]})
\cong  H^0 (X^k, L^{\boxtimes k})^{S_{k, \epsilon}} = \Lambda^k H^0 (X, L). \end{equation} To explain this isomorphism, consider the fiber diagram
\begin{center} $\xymatrix{ {\widehat{X}}_0^k \ar[r]^{\hat{f}} \ar[d]^{\hat{g}} & X_0^k \ar [d]^{g} \\ X^{[k]}_0 \ar[r]^{f} & X_0^{(k)},}$ \end{center}
where the bottom and right arrows $f$ and $g$ are the Hilbert-Chow and $S_k$-quotient morphisms respectively. The $0$ subscripts indicate that we
look everywhere at the open subschemes of zero cycles with at least $k-1$ distinct points, which is enough to identify spaces of sections since the
complements lie in codimension at least $2$. The isomorphism \eqref{identify} is obtained by pulling back $L^{[k]}$ via $\hat{g}$ and pushing forward 
by
$\hat{f}$.

In particular,
$$H^0 (X^{[n]}, L^{[n]} ) \cong H^0 (X^n,
L^{\boxtimes n})^{S_{n, \epsilon}} \cong \Lambda^n H^0 (X, L)$$ is one-dimensional, and is spanned by the divisor 
$$\theta_L = \left \{ I_V \in X^{[n]} \, \, \text{such that} \, \, h^0 (I_V \otimes L) \neq 0 \right \}.$$

Two remarks are now in order. First, under the rational map
$$\tau: X^{[k]} \times X^{[n-k]} \dashrightarrow X^{[n]}, \, \, \, (I_Z, I_W) \mapsto
I_Z \otimes I_W, $$
the rational pullback $\tau^{\star} \theta_L$ corresponds unambiguously to $\theta_{L,k}\subset X^{[k]} \times X^{[n-k]},$
and thus gives an injective map
$$\tau^{\star}: H^0 (X^{[n]}, L^{[n]}) \rightarrow H^0 (X^{[k]}, L^{[k]} ) \otimes H^0 (X^{[n-k]}, L^{[n-k]}).$$

Secondly, if $s_1, \ldots, s_n$ is a basis for $H^0 (L)$, then the unique divisor $\theta_L$ of $L^{[n]}$ corresponds up to scalars to $s_1
\wedge \cdots \wedge s_n$.  Furthermore, $I_V$ is a point in $\theta_L$ if and only if $g^{-1} (f (I_V))$ is in the vanishing locus of $s_1 \wedge
\cdots \wedge s_n$; the latter is regarded here as an element of $H^0 (X^n, L^{\boxtimes n})$ {\it{i.e.,}} is viewed as the antisymmetrization of $s_1
\otimes \cdots \otimes s_n$ in $H^0 (X^n, L^{\boxtimes n}).$ So the vanishing locus of $s_1 \wedge \cdots \wedge s_n \in H^0 (X^n, L^{\boxtimes n})$
agrees up to the quotient by the symmetric group with the vanishing locus $\theta_L$.

We denote the above inclusion map by $$\iota: H^0 (X^{[n]}, L^{[n]}) \rightarrow H^0 (X^n, L^{\boxtimes n}), \, \, \theta_L \mapsto s_1 \wedge
\cdots \wedge s_n,$$ and further let $$\lambda: H^0 (X^{[k]}, L^{[k]}) \otimes H^0 (X^{[n-k]}, L^{[n-k]}) = \Lambda^k H^0 (X, L ) \otimes
\Lambda^{n-k} H^0 (X, L) \hookrightarrow $$ $$\hookrightarrow H^0 (X, L)^{\otimes k} \otimes H^0 (X, L)^{\otimes (n-k)} = H^0 (X^n, L^{\boxtimes   
n})$$ be the tautological inclusion map on each of the two spaces in the tensor product, identifying a $k$-form with the antisymmetrization of the
corresponding tensor element. Note now that the diagram
\begin{center}
$\xymatrix{ H^0 (X^{[n]}, L^{[n]}) \ar[d]^{\tau^{\star}} \ar[r]^{\iota} & H^0 (X^n, L^{\boxtimes n}) \\
H^0 (X^{[k]}, L^{[k]}) \otimes H^0 (X^{[n-k]}, L^{[n-k]}) \ar[ur]^{\lambda} }$
\end{center}
commutes, since the vanishing loci of $\iota (\theta_L)$ and
$\lambda (\theta_{L,k})$ on $X^n$ coincide with $g^{-1} (f (\theta_L))$ on the open part of $X^n$ consisting
of tuples of distinct points.

The commutativity of the diagram implies now that under the isomorphism \eqref{identify}, the section
$\theta_{L,k} \in H^0 (X^{[k]}, L^{[k]}) \otimes H^0 (X^{[n-k]}, L^{[n-k]})$ is identified (up to scalars) with the image of
$s_1 \wedge \cdots \wedge s_n$
under the natural algebraic inclusion $$\Lambda^n H^0 (L) \hookrightarrow \Lambda^k H^0 (L) \otimes \Lambda^{n-k} H^0 (L).$$
The latter induces an isomorphism $$\Lambda^k H^0 (L)^{\vee} \rightarrow \Lambda^{n-k} H^0 (L),$$ therefore so does $\theta_{L,k}.$

\subsection {A rank $0$ example.}\label{hilbert3} After studying theta dualities for pairs of moduli spaces of rank $1$ sheaves on a surface $X$, we give 
an example when the two dual vectors both have rank $0$. Specifically we let 
$v = \left[\mathcal O_{Z}\right]$, with $Z$ a punctual
scheme of length $n$, and we let the orthogonal $K$-vector $w$ be the class of a rank $0$ sheaf on $X$ supported on a primitive 
divisor $D$, and having rank $1$ along $D$. Note that $\m_v \cong X^{[n]},$ and consider 
the morphism 
$$s:\m_{w}\to |D|,$$ sending a sheaf to its schematic support, which lies in the linear system $|D|$. 

The divisor $$\Theta=\{(Z, F) \text { such that } h^{0}(\mathcal O_{Z}\otimes F)\neq 0\} \hookrightarrow X^{[n]}\times \m_{w}$$ is the
pullback of the incidence divisor $$\Delta=\{(Z, \Sigma)\in X^{(n)}\times |D| \text { such that } Z\cap \Sigma \neq \emptyset\}$$ via the natural
morphism $$X^{[n]}\times \m_{w}\to X^{(n)}\times |D|.$$ This implies that $$\Theta_{v}=s^{\star}\mathcal O(n), \text { and } \Theta_{w}=D_{(n)}.$$
Therefore, $$H^{0}(\m_{v}, \Theta_{w})=H^{0}(X^{[n]}, D_{(n)})=H^{0}(X^{(n)}, D_{(n)})=H^{0}(X^{n}, D^{\boxtimes n})^{S_{n}}=\text
{Sym}^{n}H^{0}(D),$$ while $$H^{0}(\m_{w}, \Theta_{v})=H^{0}(\m_{w}, s^{\star}\mathcal O(n))=H^{0}(|D|, \mathcal O(n))=\text
{Sym}^{n}H^{0}(D)^{\vee}.$$ It is then clear that the two spaces of sections are naturally dual, with the duality induced by the divisor $\Theta$.

\section {Duality on $K3$ and abelian surfaces}

\subsection {Numerical evidence} We will collect evidence in favor of a strange duality theorem on surfaces $X$ with trivial canonical
bundle $K_{X}\cong \mathcal O_{X}$ {\it i.e.,} $K3$ or abelian surfaces. 

We will change the notation slightly, writing as customary $$v=\text
{ch} (E)\sqrt{\text {Todd }X}$$ for the Mukai vector of the sheaves $E$ in the moduli space $\m_{v}$. We will endow the cohomology of $X$ with the Mukai product, defined for
two Mukai vectors $v=(v_{0}, v_{2}, v_{4})$ and $w=(w_{0}, w_{2}, w_{4})$ by $$\langle v, w\rangle = \int_{X}
v_{2}w_{2}-v_{0}w_{4}-v_{4}w_{0}.$$ We will assume that the vector $v$ is {\it primitive and positive.} The latter requirement means that $v$ has
positive rank, or otherwise, in rank $0$, $c_{1}(v)$ is effective and $\langle v, v\rangle \neq 0, 4$. Moreover, we assume that the polarization $H$ is {\it generic}. This ensures that the
moduli space $\mathfrak M_{v}$ consists only of stable sheaves. (It is likely that these assumptions can be relaxed.) We will give explicit expressions for the Euler
characteristics $$\chi(\m_{v}, \Theta_{w}),$$ which will render obvious their symmetry in $v$ and $w$.

For both $K3$ and abelian surfaces, the calculation of the Euler characteristics is facilitated by the presence of a {\it holomorphic symplectic structure} 
on the
moduli spaces in question, as first established by Mukai \cite {mukai}. In fact, there are two basic examples of such holomorphic symplectic
structures which are relevant for the discussion at hand. The first is provided by the Hilbert scheme of points $X^{[n]}$ on a $K3$ surface $X$.
When $X$ is an abelian surface, a small variation is required, in order to obtain an {\it irreducible} symplectic structure. To this end, one
considers the addition map $$a:X^{[n]}\to X^{(n)}\to X, \, \, [Z]\to l_{1}z_{1}+\ldots + l_{m} z_{m},$$ where $Z$ is a punctual scheme supported on
$z_{1}, \ldots, z_{m}$, with lengths $l_{1}, \ldots, l_{m}$ respectively. The fibers of $a$ are termed generalized Kummer varieties, and are
irreducible holomorphic symplectic manifolds of dimension $2n-2$.

The relevance of these two examples resides in the following observations due to O'Grady and Yoshioka \cite {ogr} \cite{yoshioka}
\cite {ryoshioka}. First, when
$X$ is a $K3$ surface, the moduli space $\m_{v}$ is deformation equivalent to the Hilbert scheme of points $X^{[d_{v}]}$, with
$$d_{v}=\frac{1}{2}\langle v, v \rangle +1.$$ The Euler characteristics of the theta line bundles $\Theta_{w}$ on
$\m_{v}$ are deformation-invariant polynomials in the 
Beauville-Bogomolov form. These polynomials can therefore be calculated on the Hilbert scheme. 
The argument is presented in \cite {ogr2}, \cite {gny}, yielding the answer \begin{equation}\label{s1}\chi(\m_{v},
\Theta_{w})=\chi(\m_{w}, \Theta_{v})=\binom {d_{v}+d_{w}}{d_{v}}.\end {equation}

The situation is slightly more involved when $X$ is an abelian surface. In order to obtain an irreducible holomorphic symplectic structure, we need
to look at the Albanese morphism $\alpha$ of $\m_{v}$. To this end, write $${\bf R}\mathcal S:{\bf D}(X)\to {\bf D}(\widehat X)$$ for the
Fourier-Mukai transform on $X$ with respect to the normalized Poincar\'{e} sheaf $\mathcal P$ on $X\times \widehat X:$ $${\bf R}\mathcal S(x)={\bf
R}\text {pr}_{\widehat X!}(\mathcal P\otimes \text{pr}_{X}^{\star}x).$$ Following Yoshioka \cite {yoshioka}, we may define the following
'determinant' morphism $$\alpha=(\alpha^{+}, \alpha^{-}):\m_{v}\to \widehat {X} \times X$$ with $$\alpha^{+}(E)=\det E,\,\,\ \alpha^{-}(E)=\det
{\bf R}\mathcal S(E).$$ (The identification of the target of $\alpha^{+}$ with $\widehat X$ requires the translation by a fixed reference line
bundle $\Lambda$ with $c_{1}(\Lambda)=-c_{1}(v)$. The same remark applies to the morphism $\alpha^{-}$.) Yoshioka established that the fibers $K_{v}$
of the Albanese morphism $\alpha$ are {\it irreducible} holomorphic symplectic manifolds, deformation equivalent to the generalized Kummer varieties of
dimension $2d_{v}-4.$ 

It is shown in \cite {mo2} that $$\chi(K_{v}, \Theta_{w})=\frac{(d_{v}-1)^{2}}{d_{v}+d_{w}-2} \binom{d_{v}+d_{w}-2}{d_{v}-1}.$$ This formula is clearly not symmetric in $v$ and $w$, and in fact it is not expected to be so. Instead, {\it three} symmetric formulas are
obtained considering suitable variations of the moduli spaces involved. More precisely, let us write $\m^{+}_{v}$, $\m^{-}_{v}$ for the fibers
of the two morphisms $\alpha^{+}$, $\alpha^{-}$. Then,

\begin {theorem} \cite {mo3} The following three symmetries are valid \begin{equation}\label{s2}\chi(\m^{+}_{v}, \Theta_{w})=\chi(\m^{+}_{w},
\Theta_{v})=\frac{c_{1}(v\otimes w)^2}{2(d_{v}+d_{w}-2)} \binom{d_{v}+d_{w}-2}{d_{v}-1}.\end {equation} \begin{equation}\label{s3}\chi(\m^{-}_{v},
\Theta_{w})=\chi(\m^{-}_{w}, \Theta_{v})= \frac{c_{1}(\widehat v \otimes \widehat w)^2}{2(d_{v}+d_{w}-2)} \binom{d_{v}+d_{w}-2}{d_{v}-1}.\end
{equation} \begin{equation}\label{s4}\chi(K_{v}, \Theta_{w})=\chi(\m_{w}, \Theta_{v})= \frac{(d_{v}-1)^{2}}{d_{v}+d_{w}-2}
\binom{d_{v}+d_{w}-2}{d_{v}-1}.\end {equation} \end {theorem} The last equation was derived under the assumption that $c_{1}(v)$ and $c_{1}(w)$ are
proportional, which happens for instance if the Picard rank of $X$ is $1$; this assumption is likely unnecessary. In the second equation, the hats
decorating $v$ and $w$ denote the cohomological Fourier-Mukai transform.

The numerical coincidences \eqref{s1}, \eqref{s2}, \eqref{s3}, \eqref{s4} suggest strange duality phenomena on $K3$ and abelian surfaces.
However, in order to make use of the numerics provided by these equations, one needs to assume that the $\Theta$s have no higher cohomologies. 
This is true in many cases, but seems difficult to settle in general -- see \cite {mo3} for a discussion.
 
Nonetheless, in all four cases when the above numerical symmetries occur, the splitting \eqref{split} is easily established, cf. \cite {mo3}.
Moreover, the vanishing \eqref{vanish} is guaranteed if for instance $c_{1}(v\otimes w)\cdot H>0.$ This motivates the following \begin 
{conjecture}
Let $X$ be a $K3$ or abelian surface. Assume that $v$ and $w$ are primitive, positive Mukai vectors, such that $\chi(v\otimes w)=0$, and
$$c_{1}(v\otimes w)\cdot H>0.$$ Let $(\mathcal M_{v}, \mathcal M_{w})$ denote the pair $(\m_{v}, \m_{w})$ when $X$ is a $K3$-surface, or any one of
the pairs $(K_{v}, \m_{w})$, $(K_{w}, \m_{v})$, $(\m^{+}_{v}, \m^{+}_{w})$ or $(\m^{-}_{v}, \m^{-}_{w})$ when $X$ is abelian. Then, the duality
morphism $$\mathsf {D}:H^{0}(\mathcal M_{v}, \Theta_{w})^{\vee} \to H^{0}(\mathcal M_{w}, \Theta_{v})$$ is either an isomorphism or zero. \end {conjecture}

\subsection {Examples} Other than the numerical evidence provided by equations \eqref{s1}-\eqref{s4}, the conjecture has not received much 
checking. We expect that several cases can be verified, starting with the statement for Hilbert schemes established in Section 
\ref{1-1}, by applying Fourier-Mukai transformations. We present here a handful of low-dimensional examples.

\example \label{eone} {\it Intersections of quadrics.} Let us begin with a classical example, due to Mukai \cite {mukai}. Assume that $X$ is an
intersection of three smooth quadrics in $\mathbb P^{5}$, $$X=Q_{0} \cap Q_{1} \cap Q_{2}. $$ Let $C$ be a hyperplane section, represented by a
smooth curve of genus $5$. Write $\omega$ for the class of a point on $X$. Then, the moduli space of sheaves on $X$ with Mukai vector $$v=2+C+2\omega$$ is another $K3$ surface $Y$. In fact, $Y$ may be
realized as a double cover $Y\to \mathbb P^{2}$ branched along a sextic $B$ as follows. For any $\lambda=[\lambda_{0}:\lambda_{1}:\lambda_{2}] \in
\mathbb P^{2}$, let $$Q_{\lambda}=\lambda_{0}Q_{0}+\lambda_{1}Q_{1}+\lambda_{2}Q_{2}\hookrightarrow \mathbb P^{5}$$ be a quadric in the net
generated by $Q_{0}, Q_{1}, Q_{2}$. The sextic $B$ corresponds to the singular quadrics $Q_{\lambda}$. For $\lambda$ outside $B$, $Q_{\lambda}$ may
be identified with the Pl\"{u}cker embedding $$G(2, 4)\hookrightarrow \mathbb P^{5}.$$ Therefore, the tautological sequence $$0\to \mathcal A\to
\mathcal O^{\oplus 4}\to \mathcal B\to 0$$ on $G(2,4)\cong Q_{\lambda}$ gives, by restriction to $X$, two natural rank $2$ bundles $\mathcal
A^{\vee}|_{X}$ and $\mathcal B|_{X}$, both belonging to the moduli space $\m_{v}\cong Y$. The fiber of the double cover $$f:Y\to \mathbb P^{2}$$
over $\lambda\in \mathbb P^{2}\setminus B$ consists of these two sheaves on $X$.

Let $$w=1-\omega$$ be a dual vector, so that $\m_{w}\cong X^{[2]}.$ There is a dual fibration $$\widehat f: \m_{w}\to (\mathbb P^{2})^{\vee},$$
defined by assigning $$Z\in X^{[2]}\mapsto L_{Z}\in (\mathbb P^{2})^{\vee}.$$ Here $L_{Z}$ is a line in the net of quadrics, defined as $$L_{Z}=\{\lambda \in \mathbb P^{2} \text{ such that the quadric } Q_{\lambda} \text { contains the line spanned by }
Z\}.$$

Now the theta divisor $$\Theta \hookrightarrow \m_{v}\times \m_{w}$$ is the pullback of the incidence divisor
$$\Delta\hookrightarrow \mathbb P^{2}\times (\mathbb P^{2})^{\vee}$$ via the morphism $$f\times \widehat f: \m_{v}\times \m_{w}\to \mathbb
P^{2}\times (\mathbb P^{2})^{\vee}.$$ A direct argument, valid outside the branch locus, is easy to give. Any sheaf $F\in \m_{v}\setminus
f^{-1}(B)$ admits a surjective morphism $\mathcal O^{4}\to F\to 0$ which determines a morphism $X\to G(2,4)\hookrightarrow \mathbb P^{5}$. Then,
one needs to check that the statement $$h^{0}(F\otimes I_{Z})>0$$ is equivalent to the fact that the line spanned by $Z$ is contained in the quadric
$G(2,4)\hookrightarrow \mathbb P^{5}.$ A moment's thought shows this is the case, upon unravelling the definitions. We refer the reader to \cite
{ogr2}, Claim $5.16$, for a more complete argument.

This implies that $$H^{0}(\m_{v}, \Theta_{w})=H^{0}(Y, f^{\star}\mathcal O(1))\cong H^{0}(\mathbb P^{2}, \mathcal O(1))$$ is naturally dual to
$$H^{0}(\m_{w}, \Theta_{v})=H^{0}(X^{[2]}, \widehat f^{\star}\mathcal O(1))\cong H^{0}((\mathbb P^{2})^{\vee}, \mathcal O(1)).$$

\example {\it O'Grady's generalization.} Let us now explain O'Grady's generalization of this example \cite {ogr2}. This covers the case when
$$v=2+C+2\omega, \,\, w=1-\omega,$$ where now $C$ is any smooth genus $g\leq 8$ curve obtained as a hyperplane section of a {\it generic} $K3$
surface $$X\hookrightarrow \mathbb P^{g}.$$ The situation is entirely similar to what we had before, namely the two moduli spaces come equipped
with two dual fibrations $f$ and $\widehat f$.

The previous discussion goes through for the vector $w$, setting $$\widehat f:\m_{w}\cong X^{[2]}\mapsto |I_{X}(2)|^{\vee}, Z\mapsto L_{Z}$$ where
$$L_{Z}=\{\text{quadrics }Q \text { vanishing on } X, \text { and which contain the line spanned by } Z\}.$$
 
Things are more involved for the vector $v$. For generic $F\in \m_{v}$, O'Grady shows that $F$ is locally free and globally generated, and
$h^{0}(F)=4$. Therefore, there is an exact sequence $$0\to E\to H^{0}(F)\otimes \mathcal O_{X}\to F\to 0,$$ inducing a morphism $$X\rightarrow  G(2,
H^{0}(F))\cong G(2,4)\hookrightarrow \mathbb P(\Lambda^{2}H^{0}(F))\cong \mathbb P^{5}\to \mathbb P(H^{0}(\mathcal O_{X}(1)))\cong \mathbb P^{g},$$
where we used that $\Lambda^{2}F=\mathcal O_{X}(1).$ In turn, this gives a quadric $Q_{F}$ on $\mathbb P^{g}$ of rank at most $6$, vanishing on
$X$. Phrased differently, we obtain $$f:\m_{v}\mapsto |I_{X}(2)|\cong \mathbb P^{\binom{g-2}{2}-1}, \,\,\, F\mapsto Q_{F}.$$ O'Grady shows that $f$
is a morphism, which double covers its image $P$. The image $P$ is then shown to be a non-degenerate subvariety of the system of quadrics
$|I_{X}(2)|$.

As before, the theta duality is established once it is checked that the $\Theta\hookrightarrow \m_{v}\times \m_{w}$ is the pullback of the
incidence divisor via $f\times \widehat f$; see Claim $5.16$ in \cite {ogr2}.

\example {\it Isotropic Mukai vectors.} Example \ref{eone} can be generalized in a slightly different direction. The exposition below is
essentially lifted from Sawon \cite {sawon}. The idea is to exploit the situation considered in Section \ref{hilbert3}, using Fourier-Mukai to
obtain new numerics.

Assume that the Picard group of $X$ has rank $1$, and is generated by a smooth divisor $C$ with self-intersection $C^{2}=2r^{2}(g-1).$ Sawon 
studies the case
when $$v=(r, C, r(g-1)),\,\, w=(1, 0, 1-g).$$ Since the vector $v$ is isotropic, $\m_{v}$ is a new $K3$ surface $Y$. Now, $\m_{v}$ may fail to be a
fine moduli space, and therefore a universal sheaf may not exist on $X\times Y$. In fact, there is a gerbe $\alpha \in H^{2}(Y, \mathcal
O^{\star})$ which is the obstruction to the existence of a universal sheaf. In any case, there is an $\alpha$-{\it twisted} universal sheaf $\mathcal U$ on
$X\times Y$. The Fourier-Mukai transform with kernel $\mathcal U$ induces an isomorphism of moduli spaces
$$X^{[g]}\cong \m_{w}\cong \m_{\widehat w} (Y, \alpha).$$ The vector $\widehat w$ is the cohomological Fourier-Mukai dual of $w$
with kernel $\mathcal U$. Here, $\m_{\widehat w}(Y, \alpha)$ denotes the moduli space of $\alpha$-twisted sheaves on $Y$. Explicitly, each $F\in
\m_{w}$ satisfies $WIT_{1}$ with respect to $\mathcal U$, and the isomorphism is realized as $$F\mapsto \widehat F$$ where $\widehat F$ is the
non-zero cohomology sheaf of the complex ${\bf R}p_{!} (\mathcal U\otimes q^{\star}F),$ which occurs in degree $1$. Note that since $v$ and $w$ are
orthogonal, $\widehat w$ has rank $0$, and therefore there is a fibration given by taking supports $$s:\m_{w}\cong \m_{\widehat w}(Y, \alpha)\to
|D|.$$ In the above, $D$ is a smooth curve on $Y$, whose class corresponds to $w$ under the isomorphism \cite{sawon} $$H^{2}(Y)\cong v^{\perp}/v.$$ The theta 
divisor is
then realized as $$\Theta\hookrightarrow \m_{v}\times \m_{w}\cong Y\times \m_{\widehat w}(Y, \alpha),$$ as the pullback of the
incidence divisor $$\Delta\hookrightarrow Y \times |D|$$ under the natural morphism $$1\times s: Y \times \m_{\widehat w}(Y, \alpha) \to Y\times
|D|.$$ Indeed, $h^{1}(E\otimes F)\neq 0$ iff the point represented by the sheaf $[E]\in \m_{v}$ belongs to the support of the sheaf $\widehat F$.
This follows by the definition of $\widehat F$ and the base change theorem, after making use of the vanishing \eqref{vanish}.

In this case, it is clear that on $\m_{w}\cong \m_{\widehat w}(Y, \alpha)$ we have $$\Theta_{v}=s^{\star}\mathcal O(1),$$ and on $\m_{v}\cong Y$,
$$\Theta_{w}\cong D.$$ Moreover, observe that $$H^{0}(\m_{w}, \Theta_{v})=H^{0}(\m_{\widehat w}(Y, \alpha), s^{\star}\mathcal O(1))\cong H^{0}(|D|,
\mathcal O(1))=|D|^{\vee},$$ while $$H^{0}(\m_{v}, \Theta_{w})=H^{0}(Y, D)=|D|.$$ The spaces of sections are therefore naturally isomorphic, with
the isomorphism induced by the divisor $\Theta=(1 \times s)^{\star}\Delta$.

Finally, the above arguments should go through in the more general situation when $v$ is any positive isotropic vector, and
$w$ is arbitrary; the last section of \cite {sawon} contains a discussion of these matters.

\begin {example} {\it Strange duality on the projective plane.} In the light of the above discussion, one may wonder if examples of theta dualities can be established for other base
surfaces. The first obstacle in this direction is the fact that in the case of arbitrary surfaces, the required symmetry of the Euler 
characteristics has
not been proved yet, and there are no a priori reasons to expect it. We are however aware of sporadic examples of strange duality 
on ${\mathbb P}^2$, due to
D\v{a}nil\v{a} \cite {D1}\cite {D}. These examples concern the numerical classes $$\text{rk }(v)=2, c_{1}(v)=0, -19\leq \chi(v)\leq 2, \,\, 
\text {rk }(w)=0,
c_{1}(w)=1, \chi(w)=0.$$ In this case, the moduli space $\m_{w}$ is isomorphic to the dual projective space $(\mathbb P^{2})^{\vee}$, via $L\mapsto \mathcal 
O_{L}(-1)$. Moreover, $\Theta_{v}\cong \mathcal O(n)$, for $n = c_2 (v).$ For the dual moduli space, note the Barth morphism $$\mathsf J: \m_{v} \rightarrow
|\mathcal
O_{(\mathbb P^{2})^{\vee}}(n)|^{\vee},\, \, \, E \mapsto {\mathsf J}_E,$$ mapping a sheaf $E$ to its jumping set $$\mathsf J_{E}=\{\text {lines }L \text { in } 
(\mathbb
P^{2})^{\vee}\text {such that } E|_{L}\ncong \mathcal O_{L}\oplus \mathcal O_{L}\}.$$ One shows that $$\Theta_{w}=\mathsf J^{\star}\mathcal
O(1),$$ and furthermore, the duality morphism coincides with the pullback by ${\mathsf J}$, 
$$\mathsf D=\mathsf J^{\star}:H^{0}((\mathbb P^{2})^{\vee}, \mathcal O(n))^{\vee}\to H^{0}(\m_{v}, \Theta_{w}).$$
The morphism $\mathsf D$ is equivariant for the action of $SL_{3}$, and its source is an irreducible representation of the group. Therefore, to
establish that $\mathsf D$ is an isomorphism, it is enough to check the equality of dimensions for the two spaces of sections involved.

The strategy for the dimension calculation is reminiscent of Thaddeus's work on the moduli space of rank $2$ bundles over a curve. D\v{a}nil\v{a} 
considers the
moduli space of coherent systems $$0\to\mathcal O_{X}\to E\otimes \mathcal O(1),$$ where $E$ is a sheaf of rank $2$ and $c_{1}(E)=0$. Different
stability conditions indexed by a real parameter are considered; the moduli spaces undergo birational changes in codimension $2$ each time critical
values of the stability parameter are crossed. The largest critical value corresponds to a simpler space, birational to a projective bundle over a
suitable Hilbert scheme of points on $\mathbb P^{2}$. The computation can therefore be carried out on the Hilbert scheme. We refer the reader to
\cite {D1} for details. Let us finally note that further numerics $$\text{rk }(v)=2, c_{1}(v)=0, -3\leq \chi(v)\leq 2, \,\, \text {rk }(w)=0,
c_{1}(w)=2 \text { or }3, \chi (w)=0$$ were established in \cite {D}.

\end {example}

\section{Theta divisors on elliptic $K3$s}

We study here the natural theta divisor in a product of two numerically dual moduli spaces of sheaves on an elliptic $K3$ surface $X$, consisting of 
sheaves of
ranks $r$ and $s$ respectively. When $r > 2, s \geq 2,$ we show that this divisor is mapped to its counterpart on a new pair of moduli spaces, of 
sheaves with ranks $r-1$ and $s+1$,
birational to the original one via O'Grady's transformations \cite{ogr}. As the two moduli spaces are further identified birationally, using O'Grady's recipe, 
with Hilbert schemes of points $X^{[a]}$ and $X^{[b]}$ on $X$, the theta {\it bundles} on them are of the form $L^{[a]}$ respectively $L^{[b]}$ 
for a line bundle $L$ on $X$, whenever $r,s \geq 2.$  We conjecture that the theta {\it divisor} in the original product of 
moduli spaces of sheaves is accordingly identified with the theta divisor corresponding to $L$ in the product $X^{[a]} \times X^{[b]}$. This divisor is the subject of Section \ref{hilbert2}. Therefore, Proposition \ref{hilbert} and the conjecture imply that theta duality holds for many pairs on elliptic K3s.

\subsection{O'Grady's construction}

To start, we recall O'Grady's construction. Let $X$ be a smooth elliptic $K3$ surface with a section, and with
N\'{e}ron-Severi group $$NS (X) = \mathbb Z \sigma + \mathbb Z f,$$ where $\sigma$ and $f$ are the classes of the section and of the fiber
respectively. Note that $$ \sigma^2 = -2, \, \, \, f^2 = 0, \, \, \, \sigma f = 1.$$ Consider a Mukai vector $v$ with 
$$c_{1}(v)\cdot f=1,$$ {\it i.e.,}  
$$v = (r, \sigma + k f, p \omega) \in H^{2\star} (X),$$ for some $k,p \in \mathbb Z$. 
Pick a {\it suitable} polarization $H = \sigma + m f$ {\it e.g.}, assume that $m$ is sufficiently large. This choice of polarization ensures 
that $\m_{v}$ is a projective holomorphic symplectic manifold, consisting only of stable
sheaves.  We will denote by $2a$ the dimension of $\m_{v}$, $${\langle v, v \rangle} +2 = 2a.$$ 

O'Grady \cite{ogr} showed that $\m_v$ is birational to the Hilbert scheme $X^{[a]}$. We now describe this birational
isomorphism
which proceeds in steps successively modifying the rank of the sheaves in the moduli space. Since we are ultimately interested in theta divisors, we need to understand the birational transformations of the moduli space {\it {away from codimension two}}, so we will 
track the successive stages with some care.\vskip.1in

Note first that twisting with ${\mathcal O} (f)$ gives an isomorphism
$$\m_{v} \cong \m_{\tilde{v}}, \, \, \, \text{with} \, \, \, \tilde{v} = (r, \sigma + (k+r)f, (p+1) \omega).$$ This twist raises the Euler characteristic by 
$1$. We normalize $v$ by
requiring that $p = 1-r,$ and we denote the moduli space in this case
by $\m_r^a$. Thus, points in $\m_r^a$ have the Mukai vector
\begin{equation}
\label{mnorm}
v_{r,a} = (r, \sigma + (a - r (r-1))f, (1-r)\omega),
\end{equation}
and geometrically this normalization amounts to imposing that
$$\chi (E_r) = 1 \, \, \text{for} \, \, E_r \in \m_r^a.$$\vskip.1in

O'Grady shows that, as expected, the generic point $E_r$ of $\m_r^a$ has exactly one section,
\begin{equation}
\label{one}
h^0 (E_r) = 1, \, \, \text{and moreover}, \, \, h^0 (E_r (-f)) = 0.
\end{equation}
Keeping track of codimensions, we further have, importantly,
\begin{equation}
\label{two}
h^0 (E_r (-2f)) = 0 \, \, \text{for} \, \, E_r \, \, \text{outside a codimension 2 locus in} \, \, \m_r^a.
\end{equation}
Now stability forces the vanishing $h^2 (E_r (-2f)) = 0$ for all sheaves in $\m_r^a$. We conclude that
$$h^1 (E_r (-2f)) = -\chi (E_r (-2f)) = 1$$ outside a codimension 2 locus in $\m_r^a.$ O'Grady singles out an 
open subscheme $U_{r}^{a}\subset \m_{r}^{a}$ for which the vanishing \eqref{two} occurs. For sheaves $E_{r}$ in $U_{r}^{a}$ there is 
a unique nontrivial extension 
\begin{equation}
\label{basic}
0 \rightarrow {\mathcal O} \rightarrow \widetilde{E}_{r+1} \rightarrow E_{r} \otimes {\mathcal O}(-2f) \rightarrow 0.
\end{equation} 
 The resulting middle term $\widetilde{E}_{r+1}$ is torsion-free with Mukai vector $v_{r+1, a}.$ \vskip.1in

However, $\widetilde {E}_{r+1}$ may not be stable. In fact, it fails to be stable if $E_{r}$ belongs to a divisor $D_{r}$ in $U_r^a.$ For sheaves $E_{r}$ away 
from $D_{r}$, 
we set $$E_{r+1}\cong \widetilde E_{r+1}.$$ For sheaves $E_{r}$ in $D_{r}$, a stable sheaf $E_{r+1}$ is obtained by modifying $\widetilde 
E_{r+1}$.  
For $r\geq 2,$ the corresponding extension $\widetilde E_{r+1}$ has a natural rank $r$ subsheaf $G_r$ such that
$$0 \rightarrow G_r \rightarrow \widetilde E_{r+1} \rightarrow {\mathcal O} (f) \rightarrow 0.$$ The stabilization $E_{r+1}$ of 
$\widetilde E_{r+1}$ then fits in an exact sequence  
$$ 0 \rightarrow {\mathcal O} (f) \rightarrow E_{r+1} \rightarrow G_r \rightarrow 0.$$ Note that $E_{r+1}$ has Mukai vector $v_{r+1, a}$ as well. \vskip.1in

The assignment $$E_{r}\mapsto E_{r+1}$$ identifies dense open sets $U_{r}^{a}\cong U_{r+1}^{a}$ (whose complements have codimension at least
$2$) in the moduli spaces with vectors $v_{r,a}$ and $v_{r+1,a}$. This gives rise to a birational map $$\Phi_r: \m_r^a \dashrightarrow
\m_{r+1}^a.$$ The rank $1$ moduli space $\m_1^a$ is isomorphic to the Hilbert scheme $X^{[a]}$ via $$Z \mapsto I_Z (\sigma + af).$$
Additional requirements on the scheme $Z$ single out the open set $U_{1}^{a}$. For each rank, one gets therefore a birational isomorphism
of $\m_r^a$ with the Hilbert scheme $X^{[a]}.$ \vskip.1in

A  good understanding of the morphisms $\Phi_r$ hinges crucially on identifying the divisors $D_r$ along which the semistable reduction 
needs 
to be performed, and this is the most difficult part of O'Grady's work. Since the $U_r$s are isomorphic, the $D_r$s can be identified with divisors on the Hilbert scheme $X^{[a]}$. 
Let $S$ be the divisor of cycles in $X^{[a]}$ which intersect the section $\sigma$ of the elliptic fibration. In the notation of Section \ref{1-1},
$${\mathcal O} (S) = {\mathcal O}(\sigma)_{(a)}.$$ Let $T$ be the divisor consisting of points $I_Z$ such that $$h^0 (I_Z ((a-1)f)) \neq 0.$$ O'Grady proves 
that
$$D_1 = S \cup T, \, \, \text{and} \, \, D_r = S \, \, \text{for} \, \, r \geq 2.$$ 

\subsection{Theta divisors}

With these preliminaries understood, consider now two normalized moduli spaces $\m_r^a$ and $\m_s^b$ of stable sheaves on $X$, identified birationally, away from codimension 2, with Hilbert schemes. The tensor product of two points $E_r \in  
\m_r^a$ and $F_s \in \m_s^b$ has Euler characteristic
$$
\chi (E_r \otimes F_s ) = a +b - 2 - (r+s) (r+s -2).
$$
Since $$c_1 (E_r \otimes F_s). f= r+s,$$ tensorization by ${\mathcal O} (f)$ raises the Euler characteristic by $r+s$. We will assume from now on that
$$ r+s \, | \, a+ b -2.$$ In fact, we will furthermore assume that
$$ -\nu = \frac{a+b-2}{r+s} - (r+s -2)>1 \Leftrightarrow a+b\geq (r+s)^{2}+2.$$ The definition of $\nu$ is so that
$$\chi (E_r \otimes F_s \otimes {\mathcal O} (\nu f) )= 0, \, \, \text{for} \, \, E_r \in \m_r^a, F_s \in \m_s^b.$$ 

Any semistable sheaf $E$ on $X$ whose first Chern class has positive intersection with the fiber class $f$ {\it i.e,} $$c_1 (E)= \alpha \sigma + \beta
f, \text { for } \alpha > 0,$$ satisfies $H^{2}(E)=0$ forced by the stability condition. Therefore the locus \begin{equation} \label{thetaelliptic1} 
\Theta_{r,s} =
\{(E_r, F_s) \in \m_r^a \times \m_s^b \, \, \text{such that} \, \, h^0 (E_r \otimes F_s \otimes {\mathcal O} (\nu f)) \neq 0 \} \end{equation} should
indeed correspond to a divisor.  We further let $\Theta_{F_s}$ be the locus $$\Theta_{F_s} = \{E_r \in \m_{r}^a \, \, \text{such that} \, \, h^0
(E_r \otimes F_s (\nu f)) \neq 0 \}$$ in $\m_r^a$, and denote by $\Theta_{E_r}$ its analogue in $\m_s^b$.
\vskip.1in
Now Lemma I.6.19 of \cite{ogr} gives an explicit description of the morphism $$\Theta: v_{r,a}^{\perp} \rightarrow \, \text{Pic}\, (\m_r^a)\cong \text {Pic } (X^{[a]}).$$ In particular, the class of the theta line bundle ${\mathcal O} (\Theta_{F_s})$ in the Picard group of $X^{[a]}$ is \begin{equation}
\label{thetaelliptic}
{\mathcal O} (\Theta_{F_s}) = {\mathcal O} (\sigma)_{(a)}^{(r+s)} \otimes {\mathcal O} (f)_{(a)}^{2(r+s) - 2 - \nu} \otimes M, 
\, \, \text{for } \, r \geq 2, s\geq 1.
\end{equation}
As in Section \ref{hilbert1}, ${\mathcal O} (f)_{(a)}$ and ${\mathcal O} (\sigma)_{(a)}$ denote the line bundles on $X^{[a]}$ induced by the generators 
${\mathcal O} (f)$ and 
${\mathcal O} 
(\sigma)$ of the Picard group of $X$, and $-2M$ is the exceptional divisor in $X^{[a]}$. Letting 
\begin{equation}
L = {\mathcal O} ((r+s) \sigma + (2(r+s) - 2 -\nu) f) \, \, \text{on} \, \, X,
\end{equation}
equation \eqref{thetaelliptic} reads

\begin{equation}
{\mathcal O} (\Theta_{F_s}) = L^{[a]} \, \, \text{for}\, \, r\geq 2, s \geq 1.
\end{equation} 

We conclude that \begin{equation}
{\mathcal O} (\Theta_{r,s} ) = {\mathcal O} (\Theta_{r-1, s+1}) = L^{[a]} \boxtimes L^{[b]}, \, \, \text{for}  \, \, r > 2, s \geq 2.
\end{equation} The line bundle on the right comes equipped with the theta divisor $\theta_{L,a}$ discussed in Section \ref{hilbert2}. Away from 
codimension 2, this divisor is supported on the locus 
$$\{(I_Z, I_W) \in X^{[a]} \times X^{[b]} \, \, \text{such that} \, \, h^0 (I_Z \otimes I_W \otimes L) \neq 0 \}.$$
It seems now reasonable to expect that 

\begin{conjecture}
\label{thetamatch}
The locus $\Theta_{r,s}$ is a divisor. Moreover, under O'Grady's birational identification $\m_r^a \times \m_s^b  \dashrightarrow X^{[a]} \times X^{[b]},$ we have 
$$ \Theta_{r,s} = \theta_{L,a} \, \, \text{for} \, \, r,s \geq 2.$$
\end{conjecture}

Note that when $\nu<-1$, $L$ has no higher cohomology. This follows by a simple induction on $r+s$. Therefore, $\theta_{L,a}$ induces an isomorphism
$${\mathsf D}_L: H^0 (X^{[a]}, L^{[a]})^{\vee} \rightarrow H^0 (X^{[b]}, L^{[b]}).$$ With $$\Theta_{r,a} = {\mathcal O} (\Theta_{E_r}) \, \, \text{and} \, \,
\Theta_{s,b} = {\mathcal O} (\Theta_{F_s}),$$ the conjecture implies that the theta duality map 
$$H^0 (\m_r^a, \Theta_{s, b})^{\vee} \rightarrow H^0 (\m_s^b, \Theta_{r,a})$$ is an 
isomorphism. This consequence of the conjecture can be rephrased in a more intrinsic form as follows.

\begin {corollary} \label{ellipticsd} Let $v$ and $w$ be Mukai vectors of ranks $r\geq 2$ and $s\geq 2$. Assume that \begin{itemize}\item [(i)] $\chi(v\otimes w)=0$, \item [(ii)]$c_{1}(v)\cdot f=c_{1}(w)\cdot f=1$, \item [(iii)] $\langle v, v\rangle + \langle w, w\rangle\geq 2(r+s)^{2}.$\end{itemize} Then, the duality morphism $$\mathsf {D}:H^{0}(\m_{v}, \Theta_{w})^{\vee}\to H^{0}(\m_{w}, \Theta_{v})$$ is an isomorphism. 
\end{corollary}

As a first piece of evidence for Conjecture \ref{thetamatch}, we record the natural

\begin{proposition}
\label{thetaladder}
Letting $\Phi$  be the birational map $$\Phi = (\Phi_{r-1}^{-1}, \Phi_{s}) :    \m_r^a \times \m_s^b \dashrightarrow \m_{r-1}^a \times
\m_{s+1}^b,$$ and assuming that
$$ -\nu = \frac{a+b-2}{r+s} - (r+s -2) > 1,$$
we have
$$\Phi (\Theta_{r,s}) = \Theta_{r-1, s+1}, \, \, \, \text{for} \, \, \, r >2, s\geq 2.$$

\end{proposition}

{\it Proof.} It suffices to check the set-theoretic equality, since the two divisors correspond to isomorphic line bundles. The precise description of the morphisms $\Phi_{r-1}$ and $\Phi_{s}$ in the previous subsection will be crucial for establishing this fact, via a somewhat involved diagram chase. \vskip.1in

To begin, recall the
basic exact sequence \eqref{basic} of O'Grady's birational isomorphism, giving rise to two nontrivial extensions
\begin{equation}
\label{e}
0 \rightarrow {\mathcal O} \rightarrow {\widetilde{E}}_r \rightarrow E_{r-1} (-2f) \rightarrow 0
\end{equation}
and
\begin{equation}
\label{f}
0 \rightarrow {\mathcal O} \rightarrow {\widetilde{F}}_{s+1} \rightarrow F_s (-2f) \rightarrow 0.
\end{equation} Then $E_{r}$ and $F_{s+1}$ are obtained by stabilizing $\widetilde{E}_{r}$ and $\widetilde {F}_{s+1}$ if necessary; this process is required for sheaves in the divisorial locus $S$. \vskip.1in

We will first consider the situation when both $E_{r-1}$ and $F_{s}$ are chosen outside $S$, so that $$E_{r}=\widetilde E_{r}, \,\, 
F_{s+1}=\widetilde F_{s+1}.$$ Note moreover that we may assume that either $E_{r-1}$ or $F_{s}$ is locally free, as this happens outside a set 
of codimension $2$ in the product of moduli spaces.  To establish the Proposition in this case, it suffices to show that \begin{equation} 
\label{penultimate} 
h^0 (E_r \otimes F_s \otimes {\mathcal O} (\nu f)) = 0 \, \, \text{iff} \, \, h^0 (E_{r-1} \otimes F_{s+1}\otimes {\mathcal O} (\nu f)) = 0. 
\end{equation}

Tensoring \eqref{e} by $F_{s} (\nu f)$  we get the following sequence in cohomology $$ H^0 (F_{s} (\nu f)) \rightarrow H^0 (E_r \otimes F_{s} (\nu
f)) \rightarrow H^0 (E_{r-1} \otimes F_{s} ((\nu-2) f)) \stackrel{g}{\longrightarrow} H^1 (F_{s}(\nu f)) $$ $$\rightarrow H^1 (E_r \otimes F_{s}(\nu f))
\rightarrow H^1 (E_{r-1} \otimes F_{s}((\nu-2) f)) \rightarrow 0. $$ Similarly, twisting \eqref{f} by $E_{r-1} (\nu f)$ we obtain $$ H^0 (E_{r-1} (\nu f)) \rightarrow H^0 (E_{r-1} \otimes F_{s+1} (\nu
f)) \rightarrow H^0 (E_{r-1} \otimes F_{s} ((\nu-2 ) f))\stackrel{h}{\longrightarrow}  H^1 (E_{r-1}(\nu f)) $$ $$\rightarrow H^1 (E_{r-1} \otimes F_{s+1}(\nu f))
\rightarrow H^1 (E_{r-1} \otimes F_{s}((\nu -2) f)) \rightarrow 0. $$ Since $\nu < -1$, \eqref{two} implies that $$H^0 (F_{s} (\nu f)) = H^0 (E_{r-1} (\nu f)) =
0.$$ We conclude then from the two cohomology sequences that the statement \eqref{penultimate}: $$ h^0 (E_r \otimes F_s \otimes {\mathcal O} (\nu
f)) = 0 \, \, \text{iff} \, \, h^0 (E_{r-1} \otimes F_{s+1}\otimes {\mathcal O} (\nu f)) = 0$$ can be rephrased as $$g \text{ is an isomorphism iff } h \text { is an isomorphism}.$$

This last equivalence is evident when we consider the following cohomology diagram: \begin{center} $\xymatrix { H^0 (E_{r-1} \otimes F_s ((\nu -2) f) )
\ar[r]^{\hspace{.2in} h} \ar[d]^{g} & H^1 (E_{r-1} (\nu f)) \ar[d]^{\cong} \\ H^1 (F_{s} (\nu f)) \ar[r]^{\cong}& H^2 ({\mathcal O}((\nu + 2) f)) } $
\end{center} The right and bottom maps come from \eqref{e} and \eqref{f}; these morphisms are surjective. Further, the dimensions are a priori the same in all 
but the top left corner. Indeed, $$h^1 ({E}_{r-1} (\nu f)) = -\chi ({E}_{r-1} (\nu f)) = -\nu -1, \,\, h^{1}(F_{s}(\nu f)=-\chi (F_s (\nu f)) =-\nu - 1.$$ The above
equalities hold since $\chi ({E}_{r-1}) = \chi (F_s) = 1$ and twisting by $f$ raises the Euler characteristic by $1$. Also, $$h^2 ({\mathcal O} ((\nu +2)f)) =
h^0 ({\mathcal O} ((-\nu -2) f)) = -\nu -1.$$ Thus, the right vertical and bottom arrows are isomorphisms. Therefore, $g$ is an isomorphism if and only if $h$
is one. This establishes \eqref{penultimate} in the case when $E_r$ and $F_s$ arise as stable extensions.\vskip.1in

Next, we need to examine the case when $E_{r-1}$ is in the divisorial locus $S$ on $X^{[a]}$, but $F_{s+1}$ is not in the divisorial locus $S$ 
of $X^{[b]}$, so in particular \eqref{one} holds for $F_{s+1}.$ The situation when both $E_{r-1}$ and $F_{s+1}$ are in the special locus has 
codimension $2$ in the product $\m_r^a \times \m_s^b$, therefore we ignore it. For the same reason, in the arguments below we assume that 
$F_{s+1}$ is locally free.

We thus have $F_{s+1}=\widetilde F_{s+1}$. More 
delicately, to each $E_{r-1} \in D_{r-1}$, the birational map $\Phi_{r-1}$ associates a vector bundle $E_r \in D_r$ as follows. According to 
O'Grady's 
argument, 
$E_{r-1}(-2f)$ has a unique stable locally free subsheaf $G_{r-1}$ satisfying \begin{equation} \label{g} 0 \rightarrow G_{r-1} \rightarrow E_{r-1} 
(-2f) 
\rightarrow {\mathcal O}_{f_{0}} \rightarrow 0. \end{equation} In this exact sequence, the fiber $f_{0}$ is 
the unique elliptic fiber such that $$\text{dim Hom}(E_{r-1}, \mathcal O_{f_{0}})=1,$$ whereas the Hom groups with values in the structure 
sheaves of all other 
elliptic fibers are zero. Now the extension group $\text{Ext}^1 (G_{r-1}, {\mathcal O}(f))$ is two-dimensional. Among these extensions there is a unique one whose middle term has the {\it same} jumping fiber $f_0$ {\it i.e.,} we have 
 \begin{equation}\label{stabilize2} 0 \rightarrow {\mathcal O} 
(f) 
\rightarrow {E}_r \stackrel {\pi}{\rightarrow} G_{r-1} \rightarrow 0, \end{equation} and 
\begin{equation} \label{necessary} \text{dim Hom}(E_r, \mathcal O_{f_{0}})=1.\end{equation}
The emerging sheaf $E_r$ is locally free and stable. The assignment $$E_{r-1} \mapsto E_r$$ 
induces a birational 
isomorphism $D_{r-1} \dashrightarrow D_r.$

We now show that in this case also, \begin{equation} \label{ultimate0} h^0 (E_r \otimes F_s \otimes {\mathcal O} (\nu f)) = 0 \, \, 
\text{iff} \, \, h^0 (E_{r-1} \otimes F_{s+1} \otimes {\mathcal O} (\nu f)) = 0, \end{equation}
which will conclude the proof of the Proposition. \vskip.1in

Tensoring \eqref{stabilize2} by $F_s(\nu f)$ and taking cohomology, we have $$H^0 (F_s ((\nu + 1) f))\rightarrow H^0 (E_r \otimes F_s (\nu f) 
\rightarrow H^0 (G_{r-1} \otimes F_s (\nu f)) \stackrel{j}{\longrightarrow} H^1 (F_s ((\nu +1 )f))$$ $$ \rightarrow H^1 (E_r \otimes F_s (\nu f)) 
\rightarrow H^1 (G_{r-1} \otimes F_s (\nu f)) \rightarrow 0.$$ The first $H^0$ group is zero by \eqref{one}. We conclude that 
\begin{equation} \label{eqisom1} h^0 (E_r \otimes F_s (\nu f)) = 0 \, \, \text{iff } j \text { is an isomorphism}. \end{equation}

As earlier, there is a cohomology commutative diagram \begin{center} $\xymatrix {& H^0 (G_{r-1} \otimes F_s (\nu f) ) \ar[r]^{j} \ar[d]^{\beta} & H^1 (F_{s} ((\nu +1)  
f)) \ar[d]^{\cong} \\ H^{1}(E_{r}(\nu+2)f)\ar[r]^{\hspace{-.2in}\alpha}& H^1 (G_{r-1} ((\nu +2)  f)) \ar@{>>}[r] & H^2 ({\mathcal O}((\nu + 3) f)) }$ \end{center} where the left and right vertical maps are obtained from the
sequence \eqref{f}, and the top and bottom ones from \eqref{stabilize2}. The right onto vertical morphism is an isomorphism for dimension reasons. The bottom morphism is surjective. Therefore, 
\begin{equation} \label{2equiv} 
j \text { is an isomorphism iff} \, \, \overline {\beta}:H^0 (G_{r-1} \otimes F_s (\nu f)) \to \text{Coker } \alpha \text { is an isomorphism}. \end{equation} On the other hand, as in the above argument for stable extensions $E_r$ and $F_s$, \begin{equation} \label{eqisom2}
h^0 (E_{r-1} \otimes F_{s+1} (\nu f)) = 0 \, \, \text{iff} \, \, h: H^0 (E_{r-1} \otimes F_s ((\nu -2) f)) \cong H^1 (E_{r-1} (\nu f)) \end{equation} Now \eqref{eqisom1}, \eqref{2equiv}, and \eqref{eqisom2} imply assertion \eqref{ultimate0} once we establish that $$\overline {\beta} \text { is an isomorphism iff } h \text { is an isomorphism}.$$ 

This follows by chasing the commutative diagram  

\begin{center}
$\xymatrix{ &H^0 (G_{r-1} \otimes F_s (\nu f)) \ar@{^{(}->}[r] \ar[d]^{\beta} & H^0 (E_{r-1} \otimes F_s ((\nu - 2)f)) \ar[r] \ar[d]^{h}& H^0 
({F_s}\big {|}_{f_{0}} ) \ar[d]^{\cong} 
\\ H^{0}(\mathcal O_{f_{0}})\ar@{^{(}->}^{\hspace{-.3in}\gamma}[r]& H^1 (G_{r-1} ((\nu +2) f)) \ar[r] & H^1 (E_{r-1} (\nu f)) \ar@{>>}[r] & H^1 ({\mathcal O}_{f_{0}}).} $
\end{center} Here the rows are obtained from the exact sequence \eqref{g} after appropriate tensorizations, and the columns are obtained from the exact sequence 
\eqref{f}. Note that the vertical arrow on the right is an isomorphism between one-dimensional spaces. If we assume that the middle vertical arrow $h$ is an 
isomorphism, then the left vertical map $\beta$ gives an injection $$\overline {\overline {\beta}}:H^{0}(G_{r-1}\otimes F_{s}(\nu f))\to \text {Coker }\gamma,$$ and the dimension count forces it to be an isomorphism. Conversely, if $\overline {\overline {\beta}}$ is an 
isomorphism, then the dimension count shows that $$H^1 (G_{r-1} \otimes F_s (\nu f)) = 0.$$ Thus, the last map in the top row of the diagram is surjective. This implies that the middle vertical map $h$ is an 
isomorphism by the five lemma. 

To finish the proof of Proposition \ref{thetaladder}, it remains to explain that $\overline {\beta}$ and $\overline {\overline \beta}$ coincide. This 
comes down to showing that $$\text {Coker }\alpha \cong \text {Coker } \gamma.$$ We claim that there is a commutative diagram \begin {center}$\xymatrix { & 
H^{1}(E_{r}((\nu+2)f)) \ar[d]^{\alpha}\\ H^{0}(\mathcal O_{f_{0}})\ar[r]^{\hspace{-.3in}\gamma}\ar@{.>}[ur]^{\epsilon} & 
H^{1}(G_{r-1}((\nu+2)f)).}$\end {center} The map $\gamma$ is injective, and it is easy to see that the image of $\alpha$ has dimension 1. 
Once shown to exist, $\epsilon$ induces therefore an  
isomorphism between the images of $\alpha$ and $\gamma$ in $H^{1}(G_{r-1}((\nu+2)f))$.  To define $\epsilon$, we first explain that the morphism 
$$\pi:\text {Ext}^{1}(\mathcal O_{f_{0}}, E_{r}) \to \text {Ext}^{1}(\mathcal O_{f_{0}}, G_{r-1})$$ is surjective, with $\pi$ being the second map in 
\eqref{stabilize2}. Indeed, using the exact sequence \eqref{stabilize2}, we have $$\text {Ext}^{1}(\mathcal O_{f_{0}}, E_{r}) 
\stackrel{\pi}{\longrightarrow} \text {Ext}^{1}(\mathcal O_{f_{0}}, G_{r-1})\to \text{Ext}^{2}(\mathcal O_{f_{0}}, \mathcal 
O({f}))\stackrel{\tau}{\longrightarrow}\text{Ext}^{2}(\mathcal O_{f_{0}}, E_{r})\to 0.$$ It suffices to show that $\tau$ is an isomorphism. This 
follows by counting dimensions. Indeed, using equation \eqref{necessary}, we have $$\text {Ext}^{2}(\mathcal O_{f_{0}}, E_{r}) =\text 
{Ext}^{0}(E_{r}, 
\mathcal O_{f_{0}})=1, \,\,\,\text {Ext}^{2}(\mathcal O_{f_{0}}, \mathcal O(f))=1.$$
Finally, let us denote by $e\in \text {Ext}^{1}(\mathcal O_{f_{0}}, G_{r-1})$ the extension class of the exact sequence \eqref{g}. Pick $\bar e\in 
\text {Ext}^{1}(\mathcal O_{f_{0}},E_{r})$ such that $\pi(\bar e)=e$, and define $\epsilon$ to be the multiplication by $\bar e$. This choice makes 
the above triangular diagram commutative, completing the proof of the Proposition.

\end{document}